\documentclass[12pt]{amsart}
\usepackage{amssymb,latexsym}

\def\m{\mathfrak{m}}
\def\M{\mathcal{M}}
\def\P{\mathcal{P}}
\def\Pot{\rm Pot}

\newtheorem{theorem}{Theorem}

\begin{document}

\title{Mutations for quivers with potentials: Oberwolfach talk, April 2007}

\author{Andrei Zelevinsky}
\address{\noindent Department of Mathematics, Northeastern University,
  Boston, MA 02115, USA}
\email{andrei@neu.edu}

\begin{abstract}
This is an extended abstract of my talk at the
Oberwolfach Workshop ``Algebraic Groups" (April 22 - 28, 2007).
It is based on a joint work with H.~Derksen and J.~Weyman
(arXiv:0704.0649v2 [math.RA]).
\end{abstract}

\date{June 6, 2007}

\thanks{Research supported in part
by NSF (DMS) grant \# 0500534 and by a Humboldt Research Award.}

\maketitle

\makeatletter
\renewcommand{\@evenhead}{\tiny \thepage \hfill  A.~ZELEVINSKY \hfill}

\renewcommand{\@oddhead}{\tiny \hfill MUTATIONS FOR QUIVERS WITH
POTENTIALS
 \hfill \thepage}
\makeatother

This talk is based on \cite{dwz}.
We study quivers with relations given by non-commutative analogs
of Jacobian ideals in the complete path algebra. This framework
allows us to give a quiver-theoretic interpretation of
matrix mutations from the theory of cluster algebras.
This gives a far-reaching generalization of Bernstein-Gelfand-Ponomarev reflection functors.
The motivations for this work come from several sources:
superpotentials in physics, Calabi-Yau algebras, cluster algebras.

We think of  a \emph{quiver} as a quadruple $Q=(Q_0,Q_1,h,t)$ consisting
of a pair of finite sets $Q_0$ (\emph{vertices}) and $Q_1$ (\emph{arrows}), and
two maps $h:Q_1 \to Q_0$ (\emph{head}) and $t:Q_1 \to Q_0$ (\emph{tail}).
It is represented as a directed graph with the set of vertices
$Q_0$, and directed edges $a: ta \rightarrow ha$ for $a \in Q_1$.
We fix $Q_0$ once and for all, but $Q_1$ will be allowed to vary.
We also assume that $Q$ has no loops, that is, $ha \neq ta$ for
all~$a \in Q_1$.

For every $k \in Q_0$, the \emph{quiver mutation} at~$k$
transforms any quiver~$Q$ without oriented $2$-cycles
into another quiver $\overline Q = \mu_k(Q)$ of the same kind.
The quiver $\overline Q$ is obtained from $Q$ by the following three-step procedure:
\begin{enumerate}
\item For every incoming arrow $a:j \to k$ and every outgoing
arrow $b:k \to i$ in $Q$, create a ``composite" arrow $[ba]:j \to i$
in $\overline Q$.
\item Reverse all arrows at~$k$; that is, replace each arrow
$a:j \to k$ with $a^\star: k \to j$, and $b:k \to i$ with
$b^\star: i \to k$.
\item Remove any maximal disjoint union of oriented $2$-cycles.
\end{enumerate}
Note that $\mu_k$ is an involution in the following sense: $\mu_k^2(Q)$ is isomorphic to
$Q$ (non-canonically because of a choice one has to make in Step~(3)).

Our goal is to find a representation-theoretic extension of quiver mutations at arbitrary vertices.
Recall that a \emph{representation}~$M$ of~$Q$ is a family of finite-dimensional
vector spaces (over some fixed field~$K$) $(M(i))_{i \in Q_0}$
and a family of linear maps $(M(a): M(ta) \to M(ha))_{a \in Q_1}$.
If a vertex~$k$ is a \emph{sink} (i.e., there are no arrows~$b$ with
$tb = k$), or a \emph{source} (i.e., there are no arrows~$a$ with
$ha = k$), then the mutation~$\mu_k$ amounts to Step~(2).
In this case, $\mu_k$ can be extended to an operation on
representations, with the help of reflection functors.
Namely, the representation $\overline M = \mu_k(M)$ of $\overline Q$
is defined as follows.
We set $\overline M(i) = M(i)$ for $i \neq k$, and $\overline M(c) = M(c)$
for any arrow~$c$ not incident to~$k$.
If~$k$ is a sink, we set $M_{\rm in} = \bigoplus_{ha=k} M(ta)$,
assemble all maps $M(a)$ for $ha = k$ into one linear map
$\alpha: M_{\rm in} \to M(k)$, and set $\overline M(k) = \ker \alpha$.
For every $a \in Q_1$ with $ha = k$, we define the map
$\overline M(a^\star): \overline M(k) \to M(ta)$ as the embedding
$\overline M(k) \to M_{\rm in}$ followed by the projection
$M_{\rm in} \to M(ta)$.
If~$k$ is a source, we set $M_{\rm out} = \bigoplus_{tb=k} M(hb)$,
assemble all maps $M(b)$ for $tb = k$ into one linear map
$\beta: M(k) \to M_{\rm out}$, and set $\overline M(k) = {\rm coker} \ \beta$.
For every $b \in Q_1$ with $tb = k$, we define the map
$\overline M(b^\star): M(hb) \to \overline M(k)$ as the embedding
$M(hb) \to M_{\rm out}$ followed by the projection
$M_{\rm out} \to \overline M(k)$.

Note that $\mu_k^2$ is \emph{not} the identity operation.
If $k$ is a sink then $\overline {\overline M}_k = {\rm im} \
\alpha$, so in passing from $M$ to $\mu_k^2(M)$ we ``lose"
${\rm coker}\ \alpha$.
If $k$ is a source then $\overline {\overline M}_k = {\rm im} \
\beta$, so in passing from $M$ to $\mu_k^2(M)$ we ``lose" $\ker \beta$.
To remedy this, we define (following \cite{mrz}) a
\emph{decorated} representation of~$Q$ as a pair $\M = (M,V)$,
where $M$ is an ordinary representation of~$Q$, and $V = (V(i))_{i \in Q_0}$
is just a collection of finite-dimensional vector spaces attached
to the vertices (with no maps attached).
We extend the reflection functors to decorated representations by
defining $\mu_k(\M) = \overline \M = (\overline M, \overline V)$
(for $k$ a sink or a source) as follows: $\overline M(i) = M(i)$
and $\overline V(i) = V(i)$ for $i \neq k$, while
$$\overline M(k) = \ker \alpha \oplus V(k), \quad \overline V(k) =
{\rm coker}\ \alpha \quad (k \,\, \text{a sink});$$
$$\overline M(k) = {\rm coker}\ \beta \oplus V(k), \quad \overline V(k) =
\ker \beta  \quad (k \,\, \text{a source}).$$

To extend this construction to the case of an arbitrary
vertex~$k$, we encode $Q$ algebraically by the two vector
spaces $R = K^{Q_0}$ (the \emph{vertex span} of~$Q$) and $A=K^{Q_1}$ (the \emph{arrow span} of~$Q$).
The space $R$ is a commutative algebra with the $K$-basis
$\{e_i \mid i \in Q_0\}$ of minimal orthogonal idempotents adding
up to~$1$.
The space $A$ is an $R$-bimodule, with the bimodule structure
defined so that $A_{i,j} = e_i A e_j$ has the $K$-basis
$\{a \mid ha= i, ta = j\}$.
The \emph{complete path algebra} $R\langle \langle A \rangle
\rangle$ is defined as $\prod_{d\geq 0} A^{\otimes_R d}$
(with the convention $A^{\otimes_R 0} = R$);
thus, the elements of $R\langle\langle A \rangle\rangle$ are
(possibly infinite) $K$-linear combinations of \emph{paths} $a_1  \cdots a_d$ such
that all $a_k$ are arrows, and $t(a_k) = h(a_{k+1})$ for
$1 \leq k < d$.
The algebra $R\langle \langle A \rangle \rangle$ is equipped with
the $\m$-adic topology, where $\m = \m(A) = \prod_{d\geq 1} A^{\otimes_R d}$.

A \emph{potential} $S$ on $A$ is a (possibly infinite) linear
combination of \emph{cyclic} paths $a_1  \cdots a_d$ (that is,
those with $t(a_d) = h(a_1)$); since $Q$ is assumed to have no
loops, $S \in \m(A)^2$.
Let $\Pot(A)$ denote the space of potentials up to the following
\emph{cyclical equivalence}: $S \sim S'$ if
$S - S'$ lies in the closure of the $K$-span of
the elements $a_1 \cdots a_d - a_2 \cdots a_d a_1$ for all cyclic
paths $a_1 \cdots a_d$.
We call a pair $(A,S)$ with $S \in \Pot(A)$ a
\emph{quiver with potential} (QP for short).
By a \emph{right-equivalence} between QPs $(A,S)$ and $(A',S')$ we mean
an isomorphism of $R$-algebras
$\varphi: R \langle \langle A \rangle \rangle \to
R \langle \langle A' \rangle \rangle$ such that $\varphi(S)$ is
cyclically equivalent to~$S'$.

For any $K$-linear form $\xi \in A^\star$, the \emph{cyclic
derivative} $\partial_\xi$ is the continuous $K$-linear map
$\Pot(A) \to R\langle\langle A\rangle\rangle$
acting on cyclic paths by
$$\partial_\xi (a_1 \cdots a_d) =
\sum_{k=1}^d \xi(a_k) a_{k+1} \cdots a_d a_1 \cdots a_{k-1}.$$
We define the \emph{Jacobian ideal} $J(A,S)$ as the closure of
the (two-sided) ideal in $R\langle\langle A\rangle\rangle$
generated by the elements $\partial_\xi(S)$ for all $\xi \in A^\star$.
We call the quotient $R\langle\langle A\rangle\rangle/J(A,S)$
the \emph{Jacobian algebra} of~$S$, and denote
it by $\P(A,S)$.
We prove that any right-equivalence
between $(A,S)$ and $(A',S')$ sends $J(A,S)$ onto $J(A',S')$, hence induces an
isomorphism of the Jacobian algebras $\P(A,S)$ and $\P(A',S')$.

We say that a QP $(C,T)$ is \emph{trivial} if $C$ is the arrow
span of the disjoint union of oriented $2$-cycles $\{a_1, b_1\},
\dots, \{a_N, b_N\}$, and $T = b_1 a_1 + \cdots + b_N a_N$.
Then we have $J(C,T) = \m(C)$, hence $\P(C,T) = R$.

We say that a QP $(A,S)$ is \emph{reduced} if $S \in \m(A)^3$.
Our main tool in dealing with QPs and their mutations is the
following Splitting Theorem (\cite[Theorem~4.6]{dwz}.

\begin{theorem}
\label{th:trivial-reduced-splitting}
Every QP $(A,S)$ is right-equivalent to the direct sum of
a trivial QP $(A, S)_{\rm triv}$ and
a reduced QP $(A, S)_{\rm red}$.
The embedding of $(A, S)_{\rm red}$ into $(A,S)$ induces an
isomorphism of Jacobian algebras.
Furthermore, the right-equivalence classes of
$(A, S)_{\rm triv}$ and $(A, S)_{\rm red}$
are determined by the right-equivalence class of $(A,S)$.
\end{theorem}

Now we are ready to introduce the mutation of QPs at any vertex~$k$.
Assume that a QP $(A,S)$ is reduced and such that $A$ has no oriented
$2$-cycles involving~$k$; replacing~$S$ if necessary by a
cyclically equivalent potential, we can also assume that no cyclic path
occurring in~$S$ begins (and ends) at~$k$.
We define $\mu_k(A,S) = (\overline A, \overline S) = (\widetilde A, \widetilde S)_{\rm red}$,
where
\begin{itemize}
\item
$\widetilde A$ is the arrow span of the quiver obtained
from~$Q$ by the first two steps of the above three-step mutation
procedure;
\item $\widetilde S = [S] + \Delta$, where
$\Delta = \sum_{h(a) = t(b) = k} [ba]a^\star b^\star$, and $[S]$ is obtained from
$S$ by replacing each occurrence of a factor $ba$ (with $ha = tb = k$)
in a cyclic path with $[ba]$.
\end{itemize}

\begin{theorem}
\label{th:mutation-involutive}
The correspondence $\mu_k: (A,S) \to (\overline A, \overline S)$
induces an involution on the set of right-equivalence classes of reduced QPs
without oriented $2$-cycles through~$k$; that is, the
right-equivalence class of $(\overline A, \overline S)$ is
determined by that of $(A,S)$, and
$\mu_k^2(A,S)$ is right-equivalent to $(A,S)$.
\end{theorem}

Note that even if we assume that $(A,S)$ has no oriented $2$-cycles,
this may be no longer true for $(\overline A, \overline S)$.
We prove however that, for every $A$ without oriented $2$-cycles,
a generic choice of a potential~$S$ guarantees that an arbitrary
sequence of mutations can be applied to $(A,S)$ without creating
oriented $2$-cycles.

A \emph{decorated representation} of a QP $(A,S)$
is a pair $\M=(M,V)$, where $V$ is  a finite-dimensional
$R$-module, and  $M$ is a finite-dimensional $\P(A,S)$-module.
We extend the mutations of QPs to the level of their decorated
representations.
If~$A$ has no oriented cycles then $S=0$, and the decorated
representations of $(A,0)$ are just the
decorated quiver representations as defined above; furthermore, in
this case, the mutation at every sink or source coincides with the one defined above.
A right equivalence for decorated representations is defined in a natural way.
We prove that every mutation $\mu_k$ sending $(A,S)$ to $(\overline A, \overline S)$
establishes a bijection between the right-equivalence classes of
indecomposable decorated representations of $(A,S)$ and $(\overline A, \overline S)$;
furthermore, $\mu_k^2(\M)$ is right-equivalent to~$\M$.
Unfortunately, the construction of $\mu_k(\M)$ is too long to present it here.

\end{document}